\documentclass[reqno,11pt]{amsart}

\newtheorem{theorem}{Theorem}[section]
\newtheorem{lemma}[theorem]{Lemma}
\newtheorem{corollary}[theorem]{Corollary}

 \theoremstyle{remark}
\newtheorem{remark}[theorem]{Remark}

\theoremstyle{definition}

\newtheorem{assumption}[theorem]{Assumption}

\newcommand\cD{\mathcal{D}}
\newcommand\cL{\mathcal{L}}

\newcommand\bL{\mathbb{L}}
\newcommand\bH{\mathbb{H}}
\newcommand\bR{\mathbb{R}}
\newcommand\cS{\mathbb{S}}
\newcommand\ord{\text{\rm ord}\,}

\newcommand\Lie{{\rm Lie\,}}

\newcommand{\mysection}[1]{\section{#1}
      \setcounter{equation}{0}}

\begin{document}

\title[H\"ormander's theorem]{H\"ormander's theorem for parabolic equations
with  coefficients measurable in the time variable}
\author{N.V. Krylov}
\thanks{The  author was partially supported by
 NSF Grant DMS-1160569}
\email{krylov@math.umn.edu}
\address{127 Vincent Hall, University of Minnesota,
 Minneapolis, MN, 55455}

\begin{abstract}
We are dealing with possibly
degenerate second-order parabolic operators
whose coefficients are infinitely differentiable
with respect to the  space variables and only measurable
with respect to the time variable.
We impose the H\"ormander condition on the diffusion
coefficients and prove that the solutions
of the corresponding equations with right-hand sides
which are infinitely differentiable in the space
variables in a space-time domain have also this property.
\end{abstract}



\maketitle

\mysection{Introduction}

In this article we are dealing with possibly
degenerate  second-order parabolic operators
whose coefficients are infinitely differentiable
with respect to the  space variables and only measurable
with respect to the time variable. Such operators
arise, in particular, in the theory of stochastic diffusion
processes and in filtering theory
of partially observable diffusion processes.
We impose the H\"ormander condition on the diffusion
coefficients and prove that the solutions
of the corresponding equations with right-hand sides
which are infinitely differentiable in the space
variables  in a space-time domain  have also this property.
One can say that we are proving a restricted hypoellipticity
for our operators. The author intends to use this result
to prove some kind of restricted hypoellipticity
for stochastic partial differential  equations.

The problem of hypoellipticity was solved by H\"ormander
(see the references in \cite{Ho})
and attracted attention of very many researchers.
In particular, there is a probabilistic approach
to proving H\"ormander's hypoellipticity theorem
initiated by Malliavin and
extremely well presented in \cite{Ha}.

The exposition below basically follows the  lines  
 designed by H\"or\-man\-der with
substantial impact of Kohn  and Oleinik and
 Radkevic (see \cite{Ho,Ko,OR,Ra}).
It would be very interesting to find a probabilistic
proof of our results. So far, a few attempts by the author
failed   although the author of \cite{Ta} succeeded
in doing that in case the coefficients and their spatial 
derivatives are of class $C^{1}$ in $(t,x)$.
Later the probabilistic approach allowed the authors
of \cite{CaM} to weaken the continuity hypotheses with respect to $t$
to just H\"older continuity.  On the other hand, it is worth
mentioning an interesting artile \cite{CM}
where the authors do basically some of the same steps
as we do here but
under global (restricted) H\"ormander's condition and with some 
of the arguments
which the author of the present article could not quite 
follow, because some steps seem to be missing
(see, for instance, our comment 
in parentheses below \eqref{1.05.11}).
Another difference between our results and those
in \cite{CM} is that we prove infinite differentiability of
any generalized solution and not only of measure-valued ones. 

Our main result is stated in Section \ref{section 2.14.1}
and proved in Section \ref{section 1.6.1} preceded
by Section \ref{section 2.14.5} where we prove the main
a priori estimate. Section \ref{section 1.10.1}
consists of one-page collection of   well-known facts
from the theory of pseudo-differential operators.
In Section \ref{section 2.14.7} we give estimates for
the commutators of some operators and also prove a simple
and rather weak a priori estimate for parabolic 
degenerate equations.

 In conclusion we introduce some basic notation.
By $\bR^{d}$ we denote a Euclidean space of dimension $d$,
$x=(x^{1},...,x^{d})$ is a generic point of $\bR^{d}$.
All functions are assumed to be
real valued. 
  We denote
by $Du $ the gradient of $u$, $D^{2}u$ its Hessian, $D_{i}u 
=\partial u/\partial x^{i}$.
 If $a\in\bR^{d}$, we denote $L_{a}=a^{i}D_{i}$ (the summation convention
is always enforced). If $\alpha=(\alpha_{1},...,\alpha_{d})$ is a multi-index
(meaning $\alpha_{i}=0,1,...$), then
$$
D^{\alpha}:=D_{1}^{\alpha_{1}}\cdot...\cdot D_{d}^{\alpha_{d}},\quad
|\alpha|:=\alpha_{1}+...+\alpha_{d}.
$$
Finally $\partial_{t}=\partial/\partial t$.

\mysection{First steps, main ideas, and the main result}
                                      \label{section 2.14.1}

Introduce $BC^{\infty}_{b}$ as the set of real-valued or $\bR^{d}$-valued
measurable
vector-fields $\sigma$ on 
$$
Q=\{(t,x):t\in(0,1),x\in\bR^{d}\}
$$
such that for each $t\in(0,1)$,
 $\sigma (t,x)$ is infinitely differentiable with respect to $x$, and
  for any multi-index $\alpha$ we have
$$
\sup_{ t,x\in Q}|D^{\alpha}\sigma (t,x)|<\infty.
$$

Let   $d_{1}\geq1$ be an integer  and let
$$
L_{k}=\sigma^{ik}(t,x)D_{i}
=\sum_{i=1}^{d}
\sigma^{ik}(t,x)D_{i},\quad k=0,1,...,d_{1},
$$
be some given operators with  coefficients 
$\sigma^{i k}\in BC^{\infty}_{b}$.
 Define
$$
L=\partial_{t}-\sum_{k=1}^{d_{1}}L^{2}_{k}+L_{0}.
$$

 We use $(u,v)_{0}$ and $\|u\|_{0}$ for the scalar product
and the norm in $\cL_{2}=\cL_{2}(Q)$.
Set
\begin{equation}
                                                    \label{2.12.2}
H^{1,2}=\{u\in\cL_{2}:\partial_{t}u,Du,D^{2}u\in\cL_{2},
u(0+,\cdot)=0\}.
\end{equation}
To explain what we mean by $u(0+,\cdot) $ recall that
  if $u,\partial_{t}u
\in \cL_{2}$, then  there exists a $v $ such that
$u=v$ (a.e.) in $Q$ and $v(t,\cdot)$ is a continuous
$ \cL_{2} (\bR^{d})$-valued function defined on $[0,1]$.
Therefore $u(0+,\cdot)$ is well defined as $v(0 ,\cdot)$.
In the same way, $u(1-,\cdot)$ is well defined.
The following simple   fact 
is true.

\begin{lemma}
                                    \label{lemma 1.5.2}
 There is a constant $N$
such that for
 any $u\in H^{1,2}$
\begin{equation}
                                          \label{1.5.1}
\sum_{k=1}^{d_{1}}\|L_{k}u\|^{2}_{0}\leq 
(Lu,u)_{0}
+N\|u\|^{2}_{0} ,
\end{equation}
or, equivalently, for any $u\in H^{1,2}$
such that $Lu=f$
\begin{equation}
                                                           \label{1.8.2}
\sum_{k=1}^{d_{1}}\|L_{k}u\|^{2}_{0}
\leq (f,u)_{0}+N\|u\|^{2}_{0}.
\end{equation}
\end{lemma}
          
Proof. We
multiply  $Lu=f$ through by $u$ and integrate.
We get
$$
-\sum_{k=1}^{d_{1}}(L^{2}_{k}u,u)_{0}
+\tfrac{1}{2}\sum_{i=1}^{d}
\int_{Q}[\sigma^{i0}D_{i}(u^{2}) +
\partial_{t}(u^{2})]\,dxdt=(f,u)_{0}.
$$
Next, we integrate by parts and use that the derivatives of $\sigma$
are assumed to be bounded. Then we obtain
$$
|\int_{Q}\sigma^{i0}D_{i}(u^{2})\,dxdt|
\leq N\|u\|^{2}_{0},
\quad \int_{Q} 
\partial_{t}(u^{2})\,dxdt=\int_{\bR^{d}}u^{2}(1-,x)\,dx\geq0,
$$
$$
-\sum_{k=1}^{d_{1}}(L^{2}_{k}u,u)_{0}=\sum_{k=1}^{d_{1}}
(L_{k}u,L_{k}u)_{0}+\int_{\bR^{d}}uL'u
\,dxdt,
$$
where $L'u=(D_{i}\sigma^{ik})L_{k}u$. As above
$$
|\int_{Q}uL'u\,dxdt|\leq N\|u\|^{2}_{0},
$$
and to get (\ref{1.8.2}), it only remains to combine the above
results.
 The lemma is proved.

\begin{remark} 
                                           \label{remark 12.13.1}
Later on we will use the fact that the above proof 
 can be organized differently.  
We have
$$
\sum_{k=1}^{d_{1}}\|L_{k}u\|^{2}_{0}=
\sum_{k=1}^{d_{1}}(L_{k}u,L_{k}u)_{0}=(u,v)_{0},
$$
where $v=\sum_{k=1}^{d_{1}}L_{k}^{*}L_{k}u$. Obviously,
$$
v=-\sum_{k=1}^{d_{1}}L_{k}^{2}u+
\sum_{k=1}^{d_{1}}[L_{k}+L^{*}_{k}]L_{k}u=Lu-\partial_{t}u
+L'_{0}u,
$$
where $L'_{0}$ is a first-order differential operator with respect 
to $x$
 whose coefficients are in $BC^{\infty}_{b}$.  Furthermore,
$$
(L'_{0}u,u)_{0}=(u,(L'_{0})^{*}u)_{0}
=\frac{1}{2}([L'_{0}+(L'_{0})^{*}]u,u)_{0},
$$
where $[L'_{0}+(L'_{0})^{*}]$ is an operator of multiplying by 
a $BC^{\infty}_{b}$-function. Hence,
$$
|(L'_{0}u,u)_{0}|\leq N\|u\|^{2}_{0},\quad (u,v)_{0}\leq(Lu,u)_{0}+N\|u\|^{2}_{0}.
$$
\end{remark}

\begin{corollary}
                                  \label{corolary 1.5.3}
 We have
\begin{equation}
                             \label{1.6.5}
\|L_{k}u\|_{0}\leq N\|u\|^{1/2}_{0}(\|Lu\|^{1/2}_{0}+\|u\|^{1/2}_{0})
\leq N(\|Lu\|_{0}+\|u\|_{0})\quad\forall k\geq1.
\end{equation}
\end{corollary}

Indeed,  it suffices to use
(\ref{1.5.1}) and the inequality 
$$
(Lu,u)_{0}\leq\|Lu\|_{0}\,\|u\|_{0}.
$$ 

If we knew that for any $\xi\in\bR^{d}$
there exist $b^{1},...,b^{d_{1}}$ in 
$BC^{\infty}_{b}$ such that for all $i=1,...,d$
\begin{equation}
                                \label{12.30.3}
\xi^{i}=b^{k}\sigma^{ik}
\end{equation}
in $Q$, then \eqref{1.6.5} would imply
\begin{equation}
                             \label{1.6.6}
\| Du \|_{0}\leq N(\|f\|_{0}+\|u\|_{0}),
\end{equation}
where $f=Lu$. Next, one can hope to estimate second--order
derivatives of $u$ by differentiating the equation
$Lu=f$ (if $u,f$ are smooth enough),  hopefully getting a ``good'' 
equation for $u_{x}$, so that
$$
\| D^{2}u \|_{0}\leq N(\|Df \|_{0}+\|Du \|_{0})
\leq N(\|Df \|_{0}+\|f\|_{0}+\|u\|_{0}).
$$
Keeping dreaming along the same lines, one arrives at
\begin{equation}
                                         \label{2.12.1}
\|D^{\alpha} u \|_{0} \leq N(\sum_{|\beta|\leq|\alpha|}
\|D^{\beta}f \|_{0}+\|u\|_{0})
\end{equation}
for any multi-index $\alpha$.
This shows that one has a control on smoothness
of $u$ in terms of $\cL_{2}$ given that
$f$ is smooth.
 Sobolev's embedding theorems show that 
one has a control on smoothness of  $u$ in the uniform norm as well.
 This turns out to be quite sufficient for 
proving that $u$ is infinitely differentiable in $x$
with the derivatives square integrable in $t$.
Then the fact that, for each $t$, $u(t,x)$
is infinitely differentiable in $x$ follows after integrating
in $t$ the relation $Lu=f$.

It turns out that the assumption related to 
\eqref{12.30.3} can be relaxed and the estimate
as strong as \eqref{1.6.5} is not needed.
It suffices to have, 
 say $m$ $\bR^{d}$-valued 
functions $a_{1},...,a_{m}\in BC^{\infty}_{b}$ such that,
for any $\xi\in\bR^{d}$, one could find
real-valued $b_{1},...,b_{m}\in BC^{\infty}_{b}$
satisfying 
$$
\xi=b_{1}a_{1}+...+b_{m}a_{m}
$$
in $Q$ and such that 
for each $k=1,...,m$ and any 
$u\in H^{1,2}$ we have
\begin{equation}
                                \label{12.30.4}
\|L_{a_{k}}u\|_{-\delta}\leq
N(\|Lu\|_{0}+\|u\|_{0}),
\end{equation}
where $\delta\in(0,1)$ and $N$ are independent
of $k$ and $u$, and $\|\cdot\|_{-\delta}$
is the negative norm of order $-\delta$.
In that case instead of \eqref{1.6.6} we would have
$$
\| Du \|_{-\delta}\leq N(\|f\|_{0}+\|u\|_{0}),
\quad
\| u\|_{1-\delta}\leq N(\|f\|_{0}+\|u\|_{0}).
$$
By interpolation inequalities, if $\delta\leq1/2$, $\|u\|_{0}\leq
\varepsilon\|u\|_{1-\delta}+N(\varepsilon)\|u\|_{-\delta}$
for any $\varepsilon>0$. This yields
\begin{equation}
                                                 \label{4.10.1}
\| u\|_{1-\delta}\leq N(\|f\|_{0}+\|u\|_{-\delta}),
\quad
\| Du\|_{-\delta}\leq N(\|f\|_{0}+\|u\|_{-\delta}).
\end{equation}
and this 
allows for iterations as is outlined above.
Then in case $Lu=f\in BC^{\infty}_{b}$
one could again differentiate this equation
and obtain estimates of higher-order derivatives.

H\"ormander discovered that in his situation of elliptic operators
with smooth coefficients one can obtain \eqref{12.30.4}  
  for  $L_{k}$, $k=0,...,d_{1}$, for their commutators,
and then for the commutators of higher order.
Then under the condition that, a finite number
of thus obtained vector fields generates the whole space,
 one comes to global estimates like
\eqref{2.12.1} and some additional but almost standard effort
is needed in order to show that if, say, $u$ is a generalized function
in a  domain such that $Lu=0$ then $u$ is
infinitely differentiable with respect to $x$ in this domain.

As we have mentioned, we are following arguments in \cite{Ko}
and \cite{OR}.
However unlike \cite{OR}, in our setting
 we could not obtain \eqref{12.30.4} for $L_{0}$
and, therefore, we are basically bound
to the restricted version of arguments in \cite{OR}
mimicking those in \cite{Ko}. 
Accordingly, set
$\bL_{0}=\{L_{1},...,L_{d_{1}}\}$,  
$$
\bL_{n+1}= \bL_{n} \cup\{[L_{k},M]:k=1,...,d_{1},
M\in\bL_{n}\},\quad n\geq0,
$$ 
  where $[L_{k},M]=L_{k}M-ML_{k}$. Also we denote by
$\Lie_{n}$ the set of (finite) linear combinations
of elements of $\bL_{n}$ with real-valued
coefficients of class $BC^{\infty}_{b}$.
Observe that the operator $L_{0}$ is {\em not\/}
explicitly included
into $\Lie_{n}$. Fix a domain $G\subset Q$.

Everywhere in the paper we impose the following.
 
\begin{assumption}
                                             \label{assumption 1.5.1}
There exists an $n\in\{0,1,...\}$ such that
for any   $\zeta\in C^{\infty}_{0}(G)$
 we have  $\zeta D_{i}\in\Lie_{n}$ for any $i=1,...,d$. 
\end{assumption}

In order to state our main result we introduce the necessary
function spaces. Let $\cD(Q)$ be the set of generalized functions 
on $Q$. We work in the usual scale of Sobolev-Hilbert spaces
defined for any $m\in\bR$ by
$$
 H^{m}=\Lambda^{-m}\cL_{2}(\bR^{d}),
\quad\|u\|_{H^{m}}=\|\Lambda^{m}u\|_{H^{0}},
$$
where
$$
\Lambda=(1-\Delta)^{1/2}.
$$
For $m\in\bR$ define
$$
\bH^{ m}=\{u\in\cD(Q):\Lambda^{m}u \in  
\cL_{2}\},\quad
\|u\|_{m}=\|\Lambda^{m}u\|_{0}.
$$
Also introduce
$\bH^{1,m}$ as the set
of functions $u\in \bH^{m}$ such that $u(t,\cdot)
\in H^{m-1}$ for any $t\in(0,1)$ and
  there exists an $f\in \bH^{m-2}$
 such that for any $\phi\in C^{\infty}_{0}(\bR^{d})$
and any $t\in(0,1)$ we have
\begin{equation}
                                                    \label{2.14.2}
(u(t,\cdot),\phi) =\int_{0}^{t}(f(s,\cdot),\phi) 
\,ds.
\end{equation}
In case \eqref{2.14.2} holds we, naturally, write $\partial_{t}u=f$.
Briefly, one can write that
\begin{equation}
                                                 \label{4.12.2}
\bH^{1,m}=\{u\in \cL_{2}((0,1),H^{m }
 ):\partial_{t}u\in \cL_{2}((0,1),H^{m-2}
 ),u(0+,\cdot)=0\}.
\end{equation}
The above more detailed definition just makes it precise
what we mean by $\partial_{t}u$ and also emphasizes
the fact that for $u\in \bH^{1,m}$ the distributions
$u(t,\cdot)$ are uniquely defined for any $t\in(0,1)$.
Observe that what was said above about
functions in $H^{1,2}$ remains true for functions
in $\bH^{1,2}$, the latter being just
 the collection of the modifications
of elements of $H^{1,2}$ as described before 
Lemma \ref{lemma 1.5.2} or in the following remark.
\begin{remark}
                                           \label{remark 4.12.2}
One knows (see, for instance, Theorem 3 in 
\S 5.9.2 of \cite{Ev}) that if 
$u\in \cL_{2}((0,1),H^{m }
 )$ and there is an $f\in \bH^{m-2}$ such that,
for any $\phi\in C^{\infty}_{0}(\bR^{d})$, equation
\eqref{2.14.2} holds for almost all $t\in(0,1)$, then there exists
$v\in\bH^{1,m}$ such that $v(t,\cdot)$ is a uniformly
continuous $H^{m-1}$-valued function on $(0,1)$, the distributions
$v$ and $u$ coincide on $Q$  and \eqref{2.14.2} holds for  all $t$
and $\phi\in C^{\infty}_{0}(\bR^{d})$ if we replace $u$ with $v$.
In particular, this explains that the condition
$u(0+,\cdot)=0$ in
\eqref{4.12.2} makes sense.

\end{remark}

\begin{remark}
                                           \label{remark 4.11.1}
The space $\bH^{1,m}$ is a Hilbert space with squared norm
$\|u\|^{2}_{m}+\|\partial_{t}u\|^{2}_{m-2}$. One may wonder why
the $\bH^{m-2}$-norm and not a different norm of
$\partial_{t}u$ is entering the $\bH^{1,m}$-norm
of $u$. The reason is that we are going to deal
with equations $Lu=f$ and with $ \cL_{2} $-estimates
of their spacial derivatives, say of order $m$, in terms of $f$
and lower order norms of $u$. In such situation
the $ \cL_{2} $-norm of spacial 
  derivatives of $\partial_{t}u$ of order $m-2$  is obtained
from the equation itself.

By the way, also observe that almost obviously
$\Lambda^{n}\bH^{1,m}=\bH^{1,m-n}$ and $\Lambda^{n}\partial_{t}
=\partial_{t}\Lambda^{n}$ for all $m,n\in\bR$.
\end{remark}

Here is our main result.  
If $G$ is an open subset of $\bR^{d}$, then
by $C^{\infty}_{b}(G)$ we mean the set
of infinitely differentiable functions
on $G$ each of whose derivatives of any order
is bounded in $D$.

\begin{theorem}
                            \label{theorem 1.6.01}
Let $u $ be a generalized function on $G$  such that for an
$m\in\bR$ we have $u\zeta\in\bH^{1,m}$
for any $\zeta\in C^{\infty}_{0}(G)$. 
Take a $c\in BC^{\infty}_{b}$ and assume that
 $$
\zeta(L +c)u\in \bigcap_{n}\bH^{n}
$$
for any $\zeta\in C^{\infty}_{0}(G)$. Then
\begin{equation}
                                      \label{12.31.2}
\zeta u\in \bigcap_{n}\bH^{1,n}
\end{equation}
for any $\zeta\in C^{\infty}_{0}(G)$. Furthermore,
if, for some $a,b,r\in(0,1)$ and $\Gamma=(a,b)\times B_{r}$,
where $B_{r}=\{x\in\bR^{d}:
|x|<r\}$, we have $\bar\Gamma\subset G$
then for any $t\in(a,b)$ we have $u(t,\cdot)
\in C^{\infty}_{b}(B_{r})$ and for any
  multi-index $\alpha$
\begin{equation}
                                      \label{2.14.3}
\sup_{(t,x)\in\Gamma}
|D^{\alpha}u(t,x)| 
+\sup_{(t,x),(s,x)\in\Gamma}
\frac{|D^{\alpha}u(t,x)-D^{\alpha}u(s,x)|}{
|t-s|^{1/2}} <\infty.
\end{equation}
\end{theorem}

\begin{remark}
                                           \label{remark 7.14.1}
It is easy to understand from our arguments
how the left-hand side of \eqref{2.14.3} can be estimated
in terms of $(L +c)u$ and $u$. Let $\zeta\in C^{\infty}_{0}(G)$
be such that $\zeta=1$ on  an open set containing
$\bar{\Gamma}$. Then it turns out that for any $\alpha$
and $k$ such that $2(k-|\alpha|-2)>d$
the left-hand side of \eqref{2.14.3} is less than a constant
 independent of $u$ times
$$
\|\zeta(L +c)u\|_{k}+\| \zeta u\|_{m}
$$

\end{remark}

\mysection{Pseudo--differential operators}
                                       \label{section 1.10.1}

 Let $m\in\bR$ and let $A$ be a linear operator
defined on $\cup_{n}H^{n} $, mapping it into itself,
and such that, for any $n\in\bR$, it is a bounded operator mapping 
  $H ^{n}$  into $H ^{n+m}$. Then we say that $A$ is an operator
of order (at most) $m$ and write $\ord A=m$.

There is a theory of so-called pseudo--differential operators
(see \cite{Ho}). We  will be most interested in
particlar cases of such operators given by 

(i) $\Lambda^{m}$,
which is a pseudo--differential operator of order $m$,

(ii)
the pseudo--differential oprator of order zero which
is multiplication by an infinitely differentiable
 function on $\bR^{d}$, whose any derivative of any order is bounded,

(iii) the first order pseudo--differential operators $D_{i}$, $i=1,...,d$,

(iv) products of not more than seven of the above operators, and
their 
finite linear combinations.

Denote by $\cS^{m}$ the set of pseudo-differential operators of
order $m$. 
and recall   a few   facts from the theory
of pseudo--differential operators.
We borrow the next lemma from \cite{Ho}.

\begin{lemma}
                                             \label{lemma 1.5.1}
 (i). If $A\in\cS^{m}$, then $\Lambda^{-m}A$,
  $A\Lambda^{-m}\in\cS^{0}$, that is they
 are bounded operators on $H^{s}_{2}$
for any $s \in \bR$.

(ii). If $A_{1}\in\cS^{m_{1}}$ and $A_{2}\in\cS^{m_{2}}$,
then  $A_{1}A_{2}\in\cS^{m_{1}+m_{2}}$ and
 $[A_{1},A_{2}]:=A_{1}A_{2}-A_{2}A_{1}\in\cS^{m_{1}+m_{2}-1}$.
\end{lemma}

We also use a result on pointwise multipliers (see, for
instance, \cite{Tr}).
\begin{lemma}
                                                 \label{lemma 1.22.2}
Let $m>0$ and $a$ be a real-valued function of class $C^{m}_{b}
(\bR^{d})$.
Then for any $n\in(-m,m)$ there exists a constant $N$ such that for
any $u\in H^{n}$ we have
$$
\|au\|_{H^{n}}\leq N\|a\|_{C^{m}(\bR^{d})}\|u\|_{H^{n}}.
$$
\end{lemma}

\mysection{Preliminary estimates}
                                           \label{section 2.14.7}

Here is a result of simple manipulations.

\begin{lemma} 
                                                 \label{lemma 2.15.1}
 Let $a$ and $b$ be $\bR^{d}$-valued $C^{\infty}_{b}(\bR^{d})$ 
functions, $n\in\bR$, and let $A\in\cS^{n}$. 
Then there is a constant $N$ such that for any 
$u\in H^{2} \cap H^{n} $
\begin{equation}
                                                  \label{2.15.2}
|(L_{a}L_{b}u,A u)_{H^{0}}|\leq N\|L_{b}u\|_{H^{0}}
(\|L_{a}u\|_{H^{n}}+\|u\|_{H^{n}}),
\end{equation}
\begin{equation}
                                                  \label{2.15.3}
|(L_{b}L_{a}u,A u)_{H^{0}}|\leq N\|L_{a}u\|_{H^{n}}
(\|L_{b}u\|_{H^{0}}+\|u\|_{H^{0}}).
\end{equation}

\end{lemma}

Proof. We have
$$
|(L_{a}L_{b}u,A u)_{H^{0}}|=|(L_{b}u,L_{a}^{*}A u)_{H^{0}}|
\leq \|L_{b}u\|_{H^{0}}\|L_{a}^{*}A u \|_{H^{0}},
$$
where, owing to the fact that $\ord [L_{a}^{*},A ]\leq n$
and $L_{a}^{*}u=-L_{a}u+cu$ with $c\in C^{\infty}_{b}(\bR^{d})$,
$$
\|L_{a}^{*}A u \|_{H^{0}}\leq\|AL_{a}^{*} u \|_{H^{0}}
+\|[L_{a}^{*},A ]u \|_{H^{0}}\leq N
\|L _{a}u\|_{H^{n}}+N\|u\|_{H^{n}}.
$$
This proves \eqref{2.15.2}.

Next,
$$
|(L_{b}L_{a}u,A u)_{H^{0}}|=|(L_{a}u,L_{b}^{*}A u)_{H^{0}}|
\leq \|L_{a}u\|_{H^{n}}\|L_{b}^{*}A u\|_{H^{-n}},
$$
where  
$$
\|L_{b}^{*}A u\|_{H^{-n}}\leq\|AL_{b}^{*} u\|_{H^{-n}}
+\|[L_{b}^{*},A] u\|_{H^{-n}}
\leq N(\|L_{b}u\|_{H^{0}}+\|u\|_{H^{0}})
$$
and the lemma is proved.

Now comes the key estimate for
$[L_{a},L_{b}]u$.
\begin{lemma} 
                                                 \label{lemma 1.5.5}
 Let $a$ and $b$ be as in Lemma \ref{lemma 2.15.1}
and $\varepsilon\leq 1$. 
Then there is a constant $N$ such that
 for any $u\in H^{2} $
$$
\|\,[L_{a},L_{b}]u\|_{H^{\varepsilon/2-1}} 
\leq N( \|L_{a}u\| 
_{H^{\varepsilon-1}} +\|L_{b}u\| _{H^{0}}+\|u\|_{H^{0}}).
$$

\end{lemma}

Proof. We proceed as in Remark \ref{remark 12.13.1}, introduce
$$
A 
=\Lambda^{\varepsilon-2}[L_{a},L_{b}],
$$
 and observe that $\ord A \leq\varepsilon-1\leq0$ and
$$
\|\,[L_{a},L_{b}]u\|_{H^{\varepsilon/2-1}}^{2}=  
([L_{a},L_{b}]u,A u)_{H^{0}}
=(L_{a}L_{b}u,A u)_{H^{0}}-
(L_{b}L_{a}u,A u)_{H^{0}} .
$$
After that it suffices to use \eqref{2.15.2} and \eqref{2.15.3}
and the fact that $\|\cdot\|_{H^{n}}\leq\|\cdot\|_{H^{0}}$
for $n\leq0$.
The lemma is proved.

\begin{corollary}
                                                 \label{corollary 12.13.1}

If $a\in BC^{\infty}_{b}$ and for a constant $N$ 
\begin{equation}
                                                           \label{1.5.9}
\|L_{a}u\|_{ \varepsilon-1 }\leq N(\|Lu\|_{0}+\|u\|_{0})\quad\forall
u\in \bH^{1,2},
\end{equation}
then (see Corollary \ref{corolary 1.5.3})
there exists a constant $N$ such that
\begin{equation}
                                                           \label{1.5.5}
\|\,[L_{a},L_{k}]u\|_{\varepsilon/2-1}\leq N(\|Lu\|_{0}+\|u\|_{0})\quad\forall
u\in \bH^{1,2}, k=1,...,d_{1}.
\end{equation}
\end{corollary}

Since the operators $L_{k}$ satisfy (\ref{1.5.9}) (with
$\varepsilon=1$), applying repeatedly 
Corollary \ref{corollary 12.13.1}
and then using Lemma \ref{lemma 1.22.2}, we get the following.

\begin{theorem}
                                       \label{theorem 1.5.1}
Let $n\in\{0,1,...\}$ and
 $L_{a}\in\Lie_n$. Then there
are  constants $\varepsilon\in(0,1]$ and
 $N$ such that for all $u\in \bH^{1,2}$
we have
\begin{equation}
                                                  \label{12.13.1}
\|L_{a}u\|_{\varepsilon-1}\leq N(\|Lu\|_{0}+\|u\|_{0}).
\end{equation}
\end{theorem}

Estimate \eqref{12.13.1} will play the 
role of \eqref{12.30.4}
and will allow us to proceed as it is explained after
Corollary \ref{corolary 1.5.3}.

We will also use the following lemma
which does not require any H\"ormander's condition.
The lemma is quite elementary, although as happens often
with simple facts, its proof is rather long.
Before stating it we remind the reader a classical fact
(see, for instance, Section 5 of \cite{Kr99}).
\begin{theorem}
                                       \label{theorem 12.15.1}
Let $c$ be a real-valued function belonging to $BC^{\infty}_{b}$
and $\delta>0$.
Then for any $m\in\bR$ and any $f\in\bH^{m}$
there is a unique $u\in\bH^{1,m+2}$ such that
$cu+Lu+\delta \Delta u=f$.

\end{theorem}

Here is the lemma.
\begin{lemma}
                                     \label{lemma 1.7.1}
Let $c$ be a function belonging to $BC^{\infty}_{b}$
and $\delta\geq0$.
Then for any $m=0,\pm1,\pm2,...$  and $u\in \bH^{1,m+2}$
\begin{equation}
                                       \label{1.7.5}
\|u\|_{m}\leq 
N\|c u+Lu+\delta\Delta u\|_{m} ,
\end{equation}
where $N$ is independent of $u$ and $\delta$, and the set
\begin{equation}
                                       \label{1.7.6}
\{c u+Lu+\delta\Delta u:u\in  \bH^{1,m+2}\}
\end{equation}
 is everywhere dense in $\bH^{m}$.
\end{lemma} 

Proof.
First let $m\geq0$. The usual change of the unknown function
$v(t,x)=e^{\lambda t}u(t,x)$ sows that to prove \eqref{1.7.5}
it suffices to show that there are $\lambda>0$ and $N$ (independent of $u$)
such that 
\begin{equation}
                                       \label{1.7.05}
\|u\|_{m}\leq 
N\|(c+\lambda) u+Lu+\delta\Delta u\|_{m} .
\end{equation}

Take $u\in \bH^{1,m+2}$
and define $f=Lu+(c+\lambda)u$. To estimate
derivatives of order $\leq m$ of $u$ we differentiate this
equation several times and then integrate by parts.
Actually, we can make a shortcut using
Lemma \ref{lemma 1.5.2}. So,
let $\alpha$ be a multi-index with $|\alpha|\leq m$.
We have
\begin{equation}
                                        \label{1.7.3}
D^{\alpha}Lu +\lambda D^{\alpha}u+D^{\alpha}
(cu)=D^{\alpha}f.
\end{equation}
Here by usual calculus
$$
D^{\alpha}Lu+D^{\alpha}(cu)=
LD^{\alpha}u+\sum_{k\geq1}L_{k}b_{m}^{\alpha k}u+
b_{m}^{\alpha 0}u=LD^{\alpha}u+
\sum_{k\geq1}L^{*}_{k}\tilde{b}_{m}^{\alpha k}u+
\tilde{b}_{m}^{\alpha 0}u,
$$
with $b_{m}^{\alpha i}$ and $\tilde{b}_{t}^{\alpha i}$
 being certain usual differential
operators of order $\leq m$. 
Also $D^{\alpha}u\in\bH^{1,2}$. Hence from (\ref{1.7.3})
by Lemma \ref{lemma 1.5.2} we have
$$
\sum_{k\geq1}\|L_{k}D^{\alpha}u\|^{2}_{0}\leq 
(D^{\alpha}f-\lambda D^{\alpha}u -
\sum_{k\geq1}L^{*}_{k}\tilde{b}_{t}^{\alpha k}u
-\tilde{b}_{t}^{\alpha 0}u,D^{\alpha}u)_{0}+N\|u\|^{2}_{m}
$$
\begin{equation}
                                        \label{1.7.4}
\leq N\|f\|^{2}_{m}-\lambda \|D^{\alpha}u\|_{0}^{2} +
N\sum_{k\geq1}\|L_{k}D^{\alpha}u\|_{0} \|u\|_{m} 
+N\|u\|^{2}_{m}. 
\end{equation}
By remembering that $ab\leq\delta a^{2}+\delta^{-1}
b^{2}$ and using this to estimate the products
of norms in (\ref{1.7.4}), we get
$$
\lambda\|D^{\alpha}u\|_{0}^{2}\leq N\|f\|^{2}_{m}+
N\|u\|^{2}_{m}.
$$
Upon summing up with respect to $|\alpha|\leq m$, we conclude
$$
\lambda\|u\|^{2}_{m}\leq N_{1}\|f\|^{2}_{m}+
N_{1}\|u\|^{2}_{m},
$$
where $N_{1}$ is independent of $u$ and $\lambda$. 
By taking $\lambda_{0}=2N_{1}$, we finish the 
proof of (\ref{1.7.05}) and \eqref{1.7.5} for $m\geq0$
if $\delta=0$. From the above argument
 it is not hard to see that, actually,
 (\ref{1.7.05})  and \eqref{1.7.5}
 hold    for any $\delta>0$
 with the same
constants $\lambda_{0}$ and $N$ .

To prove \eqref{1.7.5} for $m\leq0$ we first prove 
the second assertion
of the lemma, which we need to do only for $\delta=0$
in light of Theorem \ref{theorem 12.15.1}.
Furthermore, since the spaces $\bH$ are nested
it suffices to prove the denseness only for $m\geq0$.
As above we need only find a $\lambda>0$ such that
\begin{equation}
                                       \label{1.7.06}
\{(c+\lambda) u+Lu:u\in  \bH^{1,m+2}\}
\end{equation}
 is everywhere dense in $\bH^{m}$.

The number $\lambda_{0}$, found above, depends on $m$, and 
we can write $\lambda_{0}=\lambda_{0}(m)$. Without loss
of generality we assume that $\lambda_{0}(m)$ is an increasing
function of $m\geq0$ and we prove
 that the set (\ref{1.7.06}) is dense in
$\bH^{m}$ for $m\geq0$ if $\lambda=\lambda_{0}(m+2)$.
 
By Theorem \ref{theorem 12.15.1}
  for $\delta>0$ 
$$
\{(c+\lambda) u+Lu+\delta \Delta u:u\in  \bH^{1,m+4}\}=\bH^{m+2}.
$$
 Therefore, for any $f\in \bH^{m+2}$ and $\delta>0$, one can find
$u^{\delta}\in  \bH^{1,m+4}$ such that
$$
 Lu^{\delta}+\delta\Delta u^{\delta}+(\lambda+c)
u^{\delta}=f.
$$
In addition,
(remember $\lambda=\lambda_{0}(m+2)$),
$$
\|u_{\delta}\|_{m+2}\leq N\|f \|_{m+2}.
$$
Hence,
$$
\|Lu^{\delta}+(\lambda+c)u^{\delta}-f\|_{m }=
 \delta\|\Delta u_{\delta}\|_{m }
\to0.
$$
This along with the fact that $\bH^{m+2}$ is dense in $\bH^{m}$
shows that (\ref{1.7.06}) is dense in $\bH^{m }$. Thus,
 \eqref{1.7.6} is  also dense in $\bH^{m}$.

Now we prove (\ref{1.7.5}) in the remaining case by using
 duality. Take $m\geq0$
and observe that for $v\in\bH^{0}$
\begin{equation}
                                                      \label{1.22.3}
\|v\|_{-m}=\sup_{\substack{f\in\bH^{m},\\
\|f\|_{m}\leq1}}(v,f)_{0}.
\end{equation}

Let $\check{\bH}^{1,m}$ be the collection of $u(1-t,x)$,
where $u\in\bH^{1,m}$. By reversing the time variable and using
the above result one easily proves that the set
$$
\{c u+L^{*}u+\delta\Delta u:u\in  \check{\bH}^{1,m+2}\}
$$
is everywhere dense in $\bH^{m}$. 
Therefore, for any $f$ with $\|f\|_{m}\leq1$ we can find
a sequence $u_{n}\in \check{\bH}^{1,m+2}$ such that
$f_{n}:=L^{*}u_{n}+c\check u_{n}\to f$ in $\bH^{m}$. By (\ref{1.7.5})
applied in reversed time
we get $\|u_{n}\|_{m}\leq N\|f_{n}\|_{m}$ with $N$ independent of
  $f$  and $u_{n}$. This proves that there is a constant $N$
such that the set $\{\|f\|_{m}\leq 1\}$ is a subset of the closure
in $\bH^{m}$ of
$$
\{L^{*}u+\delta\Delta u+cu:u\in\check{\bH}^{1,m+2},\|u\|_{m}\leq N\}.
$$
 Hence, for $v\in\bH^{1,m+2}$
$$
\|v\|_{-m}\leq\sup_{\substack{u\in \check{\bH}^{1,m+2},\\
\|u\|_{m}\leq N}}(v,
L^{*} u+\delta\Delta u+c u)_{0}
$$
$$
=\sup_{\substack{u\in \check{\bH}^{1,m+2},\\
\|u\|_{m}\leq N}}
(L v+\delta\Delta v+cv, u)_{0}\leq N\|L v+\delta\Delta v+cv\|
_{-m},
$$
and the lemma is proved.

\begin{remark}
Actually, the lemma is true for all $m$
rather than for integers only. 
However, the proof of this requires more
manipulations based on the theory of pseudo-differential
operators and is not so elementary as the above one,
the result of which is quite sufficient for our purposes. 
\end{remark}

\mysection{The main estimate in a particular case}
                                        \label{section 2.14.5}

Throughout this section we suppose
that a stronger condition than Assumption \ref{assumption 1.5.1}
is satisfied. Namely, we suppose that there exists an integer $n$
such that, $D_{i}\in\Lie_{n}$ for any $i=1,...,d$.  
  
Observe that
$$
\|u\|_{\varepsilon}^{2}=\|u\|_{\varepsilon-1}^{2}
+\sum_{i=1}^{d}\|D_{i}u\|_{\varepsilon-1}^{2}.
$$
This, together with Theorem \ref{theorem 1.5.1}, leads
to the following.
\begin{corollary}
                                             \label{corollary 1.5.6}
There
are  constants $\varepsilon\in(0,1]$ and
 $N$ such that for all $u\in \bH^{1,2}$
we have
\begin{equation}
                                                  \label{1.5.10}
\|u\|_{\varepsilon}\leq N(\|Lu\|_{0}+\|u\|_{0}).
\end{equation}
 
\end{corollary}

We thus get \eqref{12.30.4}
and we may proceed as is  explained in
Section \ref{section 2.14.1} moving to \eqref{4.10.1}
and then starting differentiating the equation in order
to obtain a priori estimates of higher order
derivatives. 

However, in Theorem \ref{theorem 1.6.01}
we are only given that $u$ is in a negative space 
and we want to show  step by step that its smoothness is
by at least $\varepsilon$ better, then by $2\varepsilon$
better and so on. That is why
we want to derive from Corollary \ref{corollary 1.5.6}
that, with the same $\varepsilon\in(0,1]$ for any $m\in\bR$,
there is a constant $N$ such that for
 all $u\in \bH^{1,m+2}$
\begin{equation}
                                                  \label{1.05.11}
\|u\|_{m+\varepsilon}\leq N(\|Lu\|_{m}+\|u\|_{m})
\end{equation}
  (this step is missing in \cite{CM} and in \cite{Ko2}).

In  
Section \ref{section 2.14.1} we explained the idea of proving (\ref{1.05.11}) 
on the basis of (\ref{1.5.10}) by differentiating the equation $Lu=f$.
Since we are interested in estimates in $H^{m}$ not only
for  $m\geq0$, we apply the operator $\Lambda^{m}$ to 
both sides of the equation $Lu=f$. Actually, this amounts to
substituting
$\Lambda^{m}u$ instead of $u$ in (\ref{1.5.10}).

If $u\in \bH^{1,m+2}$, then $\Lambda^{m}u\in \bH^{1, 2}$
 and
after substituting we get
\begin{equation}
                                                  \label{1.5.12}
\|u\|_{m+\varepsilon}\leq
 N(\|L\Lambda^{m}u\|_{0}+\|u\|_{m})
\leq N(\|L u\|_{m}+\|\,[L,\Lambda^{m}]u\|_{0}+\|u\|_{m}).
\end{equation}
Here we get into some trouble since Lemma \ref{lemma 1.5.1}
only says that $\ord [L,\Lambda^{m}]$ may be $=m+1>m+\varepsilon$,
so that we cannot absorb $\|\,[L,\Lambda^{m}]u\|_{0}$
into  either $\|u\|_{m+\varepsilon}$
 or $\|u\|_{m}$. The help comes from ``calculus'',
which shows that $[L,\Lambda^{m}]$ has a special form.  
\begin{lemma}
                                               \label{lemma 1.5.13}
If $a$ is an $\bR^{d}$-valued $C^{\infty}_{b}(\bR^{d})$ function 
and $b\in\cS^{m}$, then
$$
[L_{a}^{2},b]=b_{1}L_{a}+b_{2},
$$
where $b_{i}\in\cS^{m}$.
\end{lemma}
Indeed, $bL_{a}^{2}=L_{a}bL_{a}+cL_{a}$ and 
$L_{a}bL_{a}=L_{a}^{2}b+L_{a}c$, where $c=[b,L_{a}]$.
 Also $L_{a}c=cL_{a}+[L_{a},c]$,
where $\ord[L_{a},c]\leq\ord c\leq\ord b\leq m$.

Lemma \ref{lemma 1.5.13} allows us to organize
\eqref{1.5.12} differently.

\begin{lemma} 
                                               \label{lemma 10.28.3}
Let $m,n\in\bR$ and $A_{m}\in\cS^{m}$.
 Then
there is a constant $N$ such that
for any $u\in\bH^{1,m+n+2}$ we have
\begin{equation}
                                                     \label{10.29.1}
\|LA_{m}u\|_{n}\leq N(\|Lu\|_{m+n}+\sum_{k\geq1}
\|L_{k}u\|_{m+n}+\|u\|_{m+n}).
\end{equation}

\end{lemma}

Proof. Observe that
$$
  \|LA_{m}u\|_{n} 
\leq  \|A_{m}Lu\|_{n}+\|[A_{m},L]u\|_{n}.
$$
It follows that it only remains to estimate $[A_{m},L]u$.
However, by Lemma~\ref{lemma 1.5.13}
$$
[A_{m},L]=\sum_{k\geq1}b_{k}L_{k}+b_{0},
$$
where $\ord b_{r}=m$, $r=0,...,d_{1}$. Hence, 
$$
\|[A_{m},L]u\|_{n}
 \leq N(\sum_{k\geq1}\|L_{k}u\|_{m+n}+\|u\|_{m+n}),
$$
and the lemma is proved.

An extra term with $L_{k}u$ on the right in \eqref{10.29.1}
suggests that we  look back at (\ref{1.6.5}).
Indeed, it turns out that,  by using \eqref{1.5.1},
 one can get a somewhat stronger
estimate of $L_{k}u$ than {\em what is needed at this stage}.
We mean
\begin{equation}
                                                    \label{1.6.1}
\sum_{k\geq1}\|L_{k}u\|_{m+\varepsilon/2}
\leq  N(\|Lu\|_{m}+\|u\|_{m}),
\end{equation}
which, along with (\ref{10.29.1}) with $n=0$ and 
the first inequality in (\ref{1.5.12})
would certainly
 finish the proof of~(\ref{1.05.11}).
\begin{theorem}
                                             \label{theorem 1.6.1}
Take $\varepsilon$ from
 Corollary \ref{corollary 1.5.6}  and let $c\in 
BC^{\infty}_{b}$.
Then for any $m,n\in\bR$,
there is a constant $N$ such that for
 all $u\in \bH^{1,m+2} $
\begin{equation}
                                                  \label{1.5.11}
\|u\|_{m+\varepsilon}+
\sum_{k\geq1}\|L_{k}u\|_{m+\varepsilon/2}
\leq N(\|(L+c)u\|_{m}+\|u\|_{n}).
\end{equation}
\end{theorem}

Proof. We are going to prove that for any $m,p\in\bR$,
there is a constant $N$ such that for
 all $u\in \bH^{1,m+2} $
\begin{equation}
                                                  \label{4.10.3}
\|u\|_{m+\varepsilon}+
\sum_{k\geq1}\|L_{k}u\|_{m+\varepsilon/2}
\leq N(\|(L+c)u\|_{m}+\|u\|_{p}+
\sum_{k\geq1}\|L_{k}u\|_{p-\varepsilon/2}).
\end{equation}
This looks like a weaker estimate than \eqref{1.5.11}, but
actually by taking $p=n-1$ in \eqref{4.10.3}
and observing that 
$$
\|L_{k}u\|_{n-1-\varepsilon/2}\leq \|L_{k}u\|_{n-1}
\leq N\|u\|_{n},\quad \|u\|_{n-2}\leq\|u\|_{n}
$$
we obtain \eqref{1.5.11}.

Next, 
if we have \eqref{4.10.3} for $p=m$,
 then for larger $p$
we get it because then
$\|\cdot\|_{m}\leq\|\cdot\|_{p}$. On the other
hand, iterating \eqref{4.10.3} with $p=m$, 
we get it for all $p\leq m$.
Hence it suffices to concentrate on
 $p=m$. In this situation the observation
that $\|Lu\|_{m}\leq \|(L+c)u\|_{m}+N\|u\|_{m}$ allows us to assume that
$c\equiv0$.

We have
$$
R_{m+\varepsilon/2}:=\sum_{k\geq1}\|L_{k}u\|_{m+\varepsilon/2}\leq 
\sum_{k\geq1}\|L_{k}\Lambda^{m+\varepsilon/2}u\|_{0} 
$$
$$
+\sum_{k\geq1}\|\,[L_{k},
\Lambda^{m+\varepsilon/2}]u\|_{0}=:I_{1}+I_{2},
$$
where $\ord [L_{k},\Lambda^{m+\varepsilon/2}]\leq 
m+\varepsilon/2$,
so that by interpolation
$$
I_{2}\leq N
\|u\|_{m+\varepsilon/2}\leq N\|u\|_{m}^{1/2}
\|u\|^{1/2}_{m+\varepsilon}.
$$
Owing to (\ref{1.5.1}) and \eqref{10.29.1},  we write for $I_{1}$
$$
I_{1}^{2}\leq N\|L\Lambda^{m+\varepsilon/2}u\|_{
-\varepsilon/2}
\|\Lambda^{m+\varepsilon/2}u\|_{ \varepsilon/2}
+\|\Lambda^{m+\varepsilon/2}u\|^{2}_{0}
$$

$$
= N\|L\Lambda^{m+\varepsilon/2}u\|_{
-\varepsilon/2}\| u\|_{m+\varepsilon }
+\|u\|_{m+\varepsilon/2}^{2}
$$

$$
\leq N\| u\|_{m+\varepsilon }
(\|Lu\|_{m}+R_{m}+\|u\|_{m}).
$$

Thus,
$$
R_{m+\varepsilon/2}\leq N( \|Lu\|_{m}+R_{m}+\|u\|_{m})^{1/2}
\| u\|_{m+\varepsilon }^{1/2},
$$
which along with the first inequality in (\ref{1.5.12})
and \eqref{10.29.1} shows that
$$
\|u\|_{m+\varepsilon}+R_{m+\varepsilon/2}\leq
N(\|Lu\|_{m}+R_{m}+\|u\|_{m})
$$
$$
+N( \|Lu\|_{m}+R_{m}+\|u\|_{m})^{1/2}
\| u\|_{m+\varepsilon }^{1/2}.
$$

It follows that
$$
\|u\|_{m+\varepsilon}+R_{m+\varepsilon/2}\leq
N(\|Lu\|_{m}+R_{m}+\|u\|_{m})
$$
and since again by interpolation inequalities
$$
R_{m}\leq NR^{1/2}_{m+\varepsilon/2}R^{1/2}_{m-\varepsilon/2},
$$
we obtain \eqref{4.10.3} with $p=m$. The theorem is proved.

\mysection{Proof of Theorem \protect\ref{theorem 1.6.01}}
                                    \label{section 1.6.1}

We derive Theorem \ref{theorem 1.6.01} from  an
``interior'' version of Theorem \ref{theorem 1.6.1}.
 There is a way to localize 
the result of Theorem \ref{theorem 1.6.1}
by using methods from the theory of pseudo--differential operators
and the specific features of the problem.
We prefer to give a more universal and absolutely 
standard proof which works 
in a great variety of situations regardless
of what kind of global estimates are obtained in H\"older  
or Sobolev spaces.

We need a special cut-off function which is used
in the statement and the proof of the
 following lemma bearing on interior estimates.
Take an  infinitely differentiable function $h(p)$ of one variable
$p\in\bR$ such that $h(p)=1$ for $p\leq1$, $h(p)=0$ for
$ p \geq2$, and $0\leq h\leq1$.
For any $r>0$   define
$$
\xi_{r}(x)=h((2 |x|-r)/r),\quad \eta_{r}(t)
=h((2|t|-r)/r),\quad \zeta_{r}(t,x)
=\xi_{r}(x)\eta_{r}(t)
$$
and for $(t_{0},x_{0})\in Q$ set
$$
\zeta_{r}^{t_{0},x_{0}}(t,x)=
\zeta_{r}(t-t_{0},x-x_{0}).
$$
Observe that 
$\zeta_{r}^{t_{0},x_{0}}(t,x)=1$ 
for $(t,x)\in Q_{r}^{t_{0},x_{0}}$   and
$\zeta_{r}(t,x)=0$ 
outside $Q_{3r/2}^{t_{0},x_{0}}$, where
$$
Q_{r}^{t_{0},x_{0}}=\{(t,x):
|x-x_{0}|< r,|t-  t_{0}|< r \}.
$$

\begin{lemma}
                                   \label{lemma 1.6.3}
Let $(t_{0},x_{0})\in G$. Then
 there exist    $\varepsilon, R\in(0,1]$ such that
\begin{equation}
                                 \label{12.16.1}
Q_{6R}^{t_{0},x_{0}}\subset G
\end{equation}
 and for any
 $c\in BC^{\infty}_{b} $,
 $m,n\in\bR$, with $n\leq m$, and $r\leq R$
$$
\|\zeta^{t_{0},x_{0}}_{r}u\|_{m+\varepsilon}+
\sum_{k\geq1}\|\zeta^{t_{0},x_{0}}_{r}L_{k}u\|_{m+\varepsilon/2}
$$
\begin{equation}
                                         \label{1.6.7}
\leq Nr^{-\alpha}(\|\zeta^{t_{0},x_{0}}_{2r}(L+c)u\|_{m}+
\|\zeta^{t_{0},x_{0}}_{2r}u\|_{n})
\end{equation}
whenever $\zeta^{t_{0},x_{0}}_{2R}u\in \bH^{1,m+2}$,
where $N,\alpha>0$ are independent of $u$ and $r$
(as a matter of fact, one can take 
$\alpha=2\tau+2\tau(1+m-n)\varepsilon^{-1}$ with
 $\tau=\max(|m|,|n|)+3$).
\end{lemma}

Proof. Take $R$ so small that 
\eqref{12.16.1} holds and observe that changing $L$ outside 
$Q_{3R}^{t_{0},x_{0}}$ does not affect
$\zeta_{2r}(L+c)u$ for $r\leq R$ since $\zeta_{2r}=0$
outside $Q_{3R}^{t_{0},x_{0}}$. Bearing this in mind, take a function
$\zeta\in C^{\infty}_{0}(\bR^{d+1})$ such that $\zeta=0$
on $Q_{3R}^{t_{0},x_{0}}$,  $\zeta=1$
 outside $Q_{4R}^{t_{0},x_{0}}$, and
 $0\leq\zeta\leq1$. Define
$$
  L' _{i}= \zeta D_{i},\quad i=1,...,d,\quad
L'=\partial_{t}
-\sum_{k\geq1} L _{k} ^{2}-\sum_{i\geq1}(L' _{i})^{2}+
L_{0}.
$$
Observe that owing to Assumption
\ref{assumption 1.5.1} for any $i=1,...,d$ we have $(1-\zeta)
D_{i}\in\Lie_{n}$ and the formula $e_{i}=(1-\zeta)e_{i}
+\zeta e_{i}$ shows that  $D_{i}$ are in $\Lie_{n}$
constructed from $L_{k},L'_{j}$.   This modification of $L$
outside $Q_{3R}^{t_{0},x_{0}}$ had only one purpose to be able to
formally apply Theorem \ref{theorem 1.6.1}.

Now fix $r\leq R$ and for integers $j\geq0$ define
$$
r_{j}=r\sum_{i=0}^{j}2^{-i},\quad
\xi^{j}(x)=h(2^{j+1}(|x-x_{0}|-r_{j}+r2^{-(j+1)})/r),
$$
$$
\eta^{j}(t) =h(2^{j+1}(|t-t_{0}|-r_{j}+r2^{-(j+1)})/r),
\quad
\zeta^{j}(t,x)=\xi^{j}(x)\eta^{j}(t),
$$
so that 
$$
r_{0}=r,\quad r_{j}\uparrow2r,\quad\zeta^{0}=
\zeta_{r}^{t_{0},x_{0}}.
$$
Also, 
$\zeta^{j}\in C^{\infty}_{0}(\bR^{d+1})$,
$\zeta^{j} =1$ in $Q^{t_{0},x_{0}}_{r_{j}}$, 
$\zeta^{j} =0$ outside
$Q^{t_{0},x_{0}}_{r_{j+1}}$, and for $\tau:=\max(|m|,|n|)+3$
\begin{equation}
                                         \label{6.12.3}
\sup_{|\alpha|\leq\tau-2 ,  t,x }
|D^{\alpha}\partial_{t}\zeta^{j}|
+
\sup_{|\alpha|\leq\tau ,t,x }|D^{\alpha}\zeta^{j}|
 \leq  Nr^{-\tau}2^{\tau j},
\end{equation}
where     $N $ is independent of $j$
and $r$. All such constants below are denoted by $N$
without specifying each time that they are
independent of $j$
and $r$. Actually, \eqref{6.12.3}
holds for any $\tau$ with $N$ depending on $\tau$.
Our particular choice of it is dictated by
Lemma \ref{lemma 1.22.2}.
To finish with notation, let $f=(L+c)u$.

Since $\zeta^{t_{0},x_{0}}_{2R}u\in \bH^{1,m+2}$ by assumption,
 we can
 substitute $u\zeta^{j}$ in (\ref{1.5.11}). Then  we get
$$
\|u\zeta^{j}\|_{m+\varepsilon}+\sum_{k\geq1}\| uL_{k}\zeta^{j}
+\zeta^{j}L_{k}u \|_{m+\varepsilon/2}
$$
\begin{equation}
                                                 \label{6.12.1}
\leq
N(\|\zeta^{j}f+2\sum_{k\geq1}(L_{k}u)L_{k}\zeta^{j}
+u L \zeta^{j}\|_{m}
+\|u\zeta^{j}\|_{n}).
\end{equation}
Here owing to (\ref{6.12.3}) and  Lemma \ref{lemma 1.22.2}
$$
\|u\zeta^{j}\|_{n}=\|\zeta^{j}\{\zeta_{2r}u\}\|_{n}
\leq Nr^{-\tau}2^{\tau j}\|\zeta_{2r}u\|_{n},
$$

$$
\|\zeta^{j}f\|_{m}=\|\zeta^{j}\{\zeta_{2r}f\}\|_{m}  
\leq Nr^{-\tau}2^{\tau j}\|\zeta_{2r}f\|_{m},
$$

$$
\|u L \zeta^{j}\|_{m}=\|u\zeta^{j+1} L \zeta^{j}\|_{m}\leq
 Nr^{-\tau}2^{\tau j}\|u \zeta^{j+1}\|_{m},
$$

$$
\|(L_{k}u)L_{k}\zeta^{j}\|_{m}=
\|\zeta^{j+1}(L_{k}u)L_{k}\zeta^{j}\|_{m}\leq
Nr^{-\tau}2^{\tau j}\|\zeta^{j+1}L_{k}u \|_{m},
$$

$$
\|uL_{k}\zeta^{j}+\zeta^{j}L_{k}u\|_{m+\varepsilon/2}
\geq\|\zeta^{j}L_{k}u\|_{m+\varepsilon/2}-
Nr^{-\tau}2^{\tau j}\|u \zeta^{j+1}\|_{m+\varepsilon/2}.
$$
Hence (\ref{6.12.1}) implies that
$$
I_{j}:=\|u\zeta^{j}\|_{m+\varepsilon}+\sum_{k\geq1}\|\zeta^{j}
L_{k}u\|_{m+\varepsilon/2}
$$
\begin{equation}
                                       \label{1.7.1}
\leq N_{1}r^{-\tau}2^{\tau j}(\|\zeta_{2r}f\|_{m}+
\|\zeta_{2r}u\|_{n}+  
\|u\zeta^{j+1}\|_{m+\varepsilon/2}+\sum_{k\geq1}
\|\zeta^{j+1}L_{k}u\|_{m}).
\end{equation}
Next, we use the interpolation inequality 
$\|v\|_{k}\leq \gamma^{l-k}\|v\|_{l}+
\gamma^{p-k}\|v\|_{p}$ for any $\gamma>0$ if
 $k$ is between $l$ and $p$
(which immediately follows from the inequality
 $a^{2k}\leq a^{2l}+a^{2p}$). Then
for any $\delta>0$
$$
N_{1}\|\zeta^{j+1}L_{k}u\|_{m}\leq\delta^{\varepsilon/2}
\|\zeta^{j+1}L_{k}u\|_{m+\varepsilon/2}
+N\delta^{n-m-1}\|\zeta^{j+1}L_{k}u\|_{n-1},
$$
where again by (\ref{6.12.3})
$$
\|\zeta^{j+1}L_{k}u\|_{n-1}=\|\zeta^{j+1}L_{k}
(\zeta_{2r}u)\|_{n-1}
$$

$$
\leq Nr^{-\tau}2^{\tau j}\|L_{k}(\zeta_{2r}u)\|_{n-1}
\leq Nr^{-\tau}2^{\tau j}\|\zeta_{2r}u\|_{n}.
$$
Therefore, by introducing a parameter $\gamma>0$,
which will be specified later, 
defining $\delta$ from the 
equation $r^{-\tau}
2^{\tau j}\delta^{\varepsilon/2}=\gamma$, 
and setting
$$
\alpha=2\tau+2\tau(1+m-n)\varepsilon^{-1},
$$
we find
$$
N_{1}r^{-\tau}2^{\tau j}\|\zeta^{j+1}L_{k}u\|_{m}\leq
\gamma\|\zeta^{j+1}L_{k}u\|_{m+\varepsilon/2}
+N(\gamma)r^{-\alpha}2^{\alpha j}\|\zeta_{2r}u\|_{n},
$$
where $N(\gamma)$ depends on $\gamma$ but is independent
of $u,r,j$.

Similarly,
$$
N_{1}\|u\zeta^{j+1}\|_{m+\varepsilon/2}\leq\delta^
{\varepsilon/2}\|u\zeta^{j+1}\|_{m+\varepsilon}+
N\delta^{n-m-1}\|u\zeta^{j+1}\|_{n-1+\varepsilon/2},
$$

$$
\|u\zeta^{j+1}\|_{n-1+\varepsilon/2}\leq 
Nr^{-\tau}2^{\tau j}\|\zeta_{2r}u\|_{n-1+\varepsilon/2}\leq 
Nr^{-\tau}2^{\tau j}\|\zeta_{2r}u\|_{n},
$$

$$
N_{1}r^{-\tau}2^{\tau j}
\|u\zeta^{j+1}\|_{m+\varepsilon/2}
\leq\gamma
\|u\zeta^{j+1}\|_{m+\varepsilon}
+N(\gamma)r^{-\alpha}2^{\alpha j}\|\zeta_{2r}u\|_{n}.
$$

Hence coming back to (\ref{1.7.1}), we get 
\begin{equation}
                                          \label{12.17.1}
I_{j}\leq \gamma I_{j+1}+
N(\gamma)r^{-\alpha}2^{\alpha j}M,
\end{equation}
where $M:=\|\zeta_{2r}f\|_{m}+\|\zeta_{2r}u\|_{n}$.
Now we chose $\gamma$ so that
$$
\gamma 2^{\alpha}=1/2,
$$
multiply both parts of \eqref{12.17.1} by 
$\gamma^{j}$, sum up for
$j=0,1,2,...$. Then we obtain
\begin{equation}
                                          \label{1.7.2}
I_{0}+S\leq S+Nr^{-\alpha}M,
\end{equation}
where
$$
S:=\sum_{j=1}^{\infty}\gamma^{j}\big(
\|u\zeta^{j}\|_{m+\varepsilon}+\sum_{k\geq1}
\|\zeta^{j}L_{k}u\|_{m+\varepsilon/2}\big),
$$
and as above in light of  (\ref{6.12.3})
$$
\|u\zeta^{j}\|_{m+\varepsilon}+\sum_{k\geq1}
\|\zeta^{j}L_{k}u\|_{m+\varepsilon/2}
\leq Nr^{-\tau}2^{\tau j}\big(
\|u\zeta_{2r}\|_{m+\varepsilon}+\sum_{k\geq1}
\|\zeta_{2r}L_{k}u\|_{m+\varepsilon/2}\big),
$$
which along with the inequality $\gamma
2^{\tau}<1$  and the fact that 
$\zeta^{t_{0},x_{0}}_{2R}u\in \bH^{ m+2}$ yield  that $S<\infty$.
This shows that (\ref{1.7.2}) coincides with  
(\ref{1.6.7}) and
the lemma is proved.

{\bf Proof of Theorem \ref{theorem 1.6.01}}. To prove the first assertion
it suffices to show that for any $(t_{0},x_{0})\in G$
there is a function $\zeta\in C^{\infty}_{0}
(G)$ which equals one in a neighborhood of $(t_{0},x_{0})$
and is such that $\zeta u\in
\bH^{1,k}$ for all $k$.

Fix a point  $(t_{0},x_{0})\in G$, take $R$ and $\varepsilon$
from Lemma \ref{lemma 1.6.3}, and  reduce $R$ if necessary so that
$$
v:=\zeta^{t_{0},x_{0}}_{R}u\in \bH^{1,m}.
$$

Next, define
$$
f=Lv+c v
$$
 and observe that,
since $v\in \bH^{1,m}$, 
we have $f\in \bH^{m-2}$. Therefore 
by Theorem \ref{theorem 12.15.1},
for $\delta>0$, there exists a unique solution
$v_{\delta}\in \bH^{1,m}$  of the equation
$$
Lv_{\delta}+\delta\Delta v_{\delta}+cv_{\delta}=f.
$$

By Lemma \ref{lemma 1.7.1}  
\begin{equation}
                                                      \label{1.10.3}
\sup_{\delta}\|v_{\delta}\|_{m-2}<\infty,
\end{equation}
\begin{equation}
                                                       \label{1.10.2}
\|v-v_{\delta}\|_{m-4}\leq N\|L(v-v_{\delta})+c
(v-v_{\delta})\|_{m-4}=
\|\delta\Delta v_{\delta}\|_{m -4}
\to0
\end{equation}
as $\delta\downarrow0$. 

Next, since on $Q^{t_{0},x_{0}}_{R}$ we have $u=v$ and $Lu=Lv$, 
it holds that
$f=Lu+cu$ on $Q^{t_{0},x_{0}}_{R}$, where by assumption 
the right-hand side on $Q^{t_{0},x_{0}}_{R}$
 is 
a restriction of an $\cap_{k}\bH^{k}$-function.
By parabolic interior regularity theory
for uniformly nondegenerate equations
(see, for instance, the proof of 
Corollary 4.2.1  of \cite{Kr08}), 
$\zeta_{r}^{t_{0},x_{0}}v_{\delta}\in \cap_{k}\bH^{1,k}$ for any $r\in(0,R)$,
which along with Lemma \ref{lemma 1.6.3} implies that,
for any $k$,
$$
\|\zeta^{t_{0},x_{0}}_{R/4}v_{\delta}\|_{k+\varepsilon}\leq
N\|\zeta^{t_{0},x_{0}}_{R/2}(Lu+cu)\|_{k}+N(1+\delta)
\|v_{\delta}\|_{m -2},
$$
where $N$ are independent of $\delta$.
Hence, by (\ref{1.10.3}) (see also Remark
\ref{remark 4.11.1}), we get that 
$\zeta^{t_{0},x_{0}}_{R/4}v_{\delta}$ are uniformly 
bounded in $\bH^{1,k}$ for $\delta\in(0,1]$
and, by (\ref{1.10.2}), that
$\zeta^{t_{0},x_{0}}_{R/4}v=\zeta^{t_{0},x_{0}}_{R/4}u\in 
\bH^{1,k}$. 
This proves the first assertion of the theorem.

The second assertion of the theorem 
follows from the first one by embedding
theorems. Indeed, the fact that
 $\zeta u\in\bH^{1,n}$ for all $n\geq1$  implies that $\partial_{t}
(\zeta u)\in\bH^{ n}$ for all $n\geq1$  and then equation
\eqref{2.14.2} implies that $\zeta(t,\cdot)u(t,\cdot)\in H^{n}$
for any $t \in(0,1)$ and   is $1/2$-H\"older
continuous with respect to $t$ in the $H^{n}$-norm. 
Since this holds for any $n$, an application
of the Sobolev embedding theorem 
with an appropriate $\zeta$ yields \eqref{2.14.3}.


\begin{thebibliography}{m}

 
\bibitem{CM} M. Chaleyat-Maurel and M. Michel,
{\em Hypoellipticity theorems and conditional laws\/}, Z. Wahrsch. 
Verw. Gebiete, Vol. 65 (1984), No. 4,
573Ð597.

\bibitem{CaM} P. Cattiaux and L. Mesnager, {\em
Hypoelliptic non-homogeneous diffusions\/},
Probab. Theory Related Fields. Vol. 123 (2002), No. 4, 453--483.

\bibitem{Ev} L.C. Evans, ``Partial Differential Equations",
Graduate Studies in Mathematics, Vol. 19,
American Mathematical Society, Providence, RI, 1998.

\bibitem{Ha} M. Hairer, {\em 
On Malliavin's proof of H\"ormander's theorem\/},
Bull. Sci. math., Vol. 135 (2011), 650--666.

\bibitem{Ho} L. H\"ormander, ``The analysis of linear
 partial differential
operators", Vol. 3, Springer, 1985. 

 
\bibitem{Ko2} J.J. Kohn, {\em Pseudo-differential operators and
hypoellipticity\/}, Proc. Sympos. Pure Math. 
(Univ. California, Berkeley, Calif., 1971),
Partial differential equations,
Vol. 23 (1973),
 61--69, Amer.
Math. Soc., Providence, R.I.  

\bibitem{Ko} J.J. Kohn, {\em Pseudo-differential operators and
non--elliptic problems\/}, pp. 157--165 in
 Pseudo-Diff. Operators (C.I.M.E., Stresa,
1968)  Edizioni Cremonese, Rome, 1969.

\bibitem{Kr99} N.V. Krylov, {\em An analytic approach to SPDEs}, 
pp. 185-242 in
Stochastic Partial Differential Equations: Six Perspectives,
Mathematical Surveys and Monographs, Vol. 64,
AMS, Providence, RI, 1999.

\bibitem{Kr08} N.V. Krylov,
``Lectures on elliptic and parabolic equations
in Sobolev spaces", Amer.
Math. Soc., Providence, RI, 2008.

\bibitem{OR} O.A. Oleinik, E.V.  Radkevich,  ``Second order 
equations 
with nonnegative characteristic form", Itogi Nauki, Mat. Analis  1969, 
VINITI, Moscow,  1971 
in Russian; English translation: Amer. Math. Soc., Plenum Press, 
Providence R.I., 1973.

\bibitem{Ra}  E.V.  Radkevich, {\em On a theorem of L. H\"ormander},
Uspekhi Matem. Nauk, Vol.~24 (1969), No.~2,
233-234 (in Russian).

 
\bibitem{Ta}  S. Taniguchi, {\em
Applications of Malliavin calculus to
time-dependent system of heat equations\/},
 Osaka J. Math., 89 (1991), 457Ð485.

 

\bibitem{Tr} H.~Triebel,  Theory of function spaces II,
Birkh\"auser Verlag, Basel--Boston--Berlin, 1992.

\end{thebibliography}
\end{document}